\documentclass[12pt]{amsart}

\usepackage{amsmath}
\usepackage{amssymb}
\usepackage{latexsym}
\usepackage{amscd}
\usepackage{graphics}

\newtheorem{Th}{Theorem}[section]
\newtheorem{lemma}[Th]{Lemma}
\newtheorem{corollary}[Th]{Corollary}
\newtheorem{example}[Th]{Example}

\newcommand{\bd}{\partial}

\newcommand{\F}{\mathbb{F}}
\newcommand{\Z}{\mathbb{Z}}
\newcommand{\Q}{\mathbb{Q}}
\newcommand{\rk}{\mathrm{rk}\,}
\newcommand{\img}{\mathrm{Im}\ }
\newcommand{\Ker}{\mathrm{Ker}\ }
\newcommand{\Lin}{\mathrm{Lin}\ }
\newcommand{\R}{\mathbb{R}}
\newcommand{\lk}{\mathrm{lk}\ }
\newcommand{\Tor}{\mathrm{Tors}\ }

\begin{document}

\title{
Embedding
$3$-manifolds with boundary\\
into closed
$3$-manifolds}

\author{Dmitry Tonkonog
}
\thanks{The author was partially supported by Euler-Program at DAAD (German Academic
Exchange Service) in 2009/2010 and by Dobrushin Scholarship at the Independent
University of Moscow in 2011}

\address{Department of Differential Geometry and Applications, 
Faculty of Mechanics and Mathematics,
Moscow State University, Moscow, 119991 Russia.}
\email{dtonkonog@gmail.com}

\begin{abstract} 

We prove  
that there is an algorithm which
determines whether or not a given 2-polyhedron can be embedded into
some integral homology 3-sphere.

This is a corollary of the following main result.
Let $M$ be a compact connected orientable
3-manifold with boundary.
Denote $G=\Z$, $G=\Z/p\Z$ or $G=\Q$.
If $H_1(M;G)\cong G^k$ and $\bd M$ is a surface of genus $g$, 
then the minimal
group $H_1(Q;G)$ for
closed 3-manifolds $Q$ containing $M$
is isomorphic to $G^{k-g}$.

Another corollary is that for a graph $L$ the minimal number
$\rk H_1(Q;\Z)$
for closed orientable 3-manifolds $Q$ containing $L\times S^1$ is twice the
orientable genus of the graph.
\end{abstract}

\begin{flushright}
To appear in Topology and its Applications (2011).
\end{flushright}
\smallskip

\maketitle

\section{Introduction and main results}

\smallskip
Let $M$ be a compact orientable 3-manifold with boundary.
In Theorem~\ref{Th1}, we
find the minimal rank of $H_1(Q;\F)$
for all closed 3-manifolds $Q$ containing $M$
(i.e. such that $M$ embeds into $Q$)
in terms of homology of $M$.
Here $\F$ is one of the fields $\Q$ or $\Z_p:=\Z/p\Z$.
Theorem \ref{Th2} is an integral version of Theorem~\ref{Th1}.
The following are two corollaries.

\begin{corollary}
\label{CorHomS}
Suppose $G=\Z$, $\Z_p$ or $\Q$. 
There exists an algorithm
that for any given (finite) 2-polyhedron $P$
tells if $P$ is embeddable into some
$G$-homology 3-sphere
(the sphere is not fixed in advance).
\end{corollary}

According to \cite{MTW11}, the existence
of an algorithm recognizing embeddability of 2-polyhedra in
$\R^3$ is unknown, cf. \cite{GS06}.

\begin{corollary}
\label{CorGrGen}
Let $L$ be a connected graph of genus $g(L)$. Suppose $\F=\Z_p$ or $\F=\Q$.
The minimal number $\dim H_1(Q;\F)$ for closed
orientable
3-manifolds $Q$ containing $L \times S^1$
equals to $2g(L)$.
\end{corollary}

Here the {\it genus of graph} $g(L)$ is the minimal $g$ such that $L$ embeds into 
a surface of genus $g$ \cite{MT01}.
To prove these corollaries, we use the classification of 
3-thickenings of 2-polyhedra \cite{BRS99, La00, Sk95}.
In particular, from the cited papers we derive Lemma~\ref{LemThick}
stating that
all orientable 3-thickenings of a given 2-polyhedron are
algorithmically constructible.

\begin{Th}
\label{Th1}
Let $M$ be a compact
connected 3-manifold with orientable boundary. Denote
$g:=\rk H_1(\bd M;\Z)/2$.
Take a field $\F=\Z_p$ or $\F=\Q$.
Suppose $M$ is orientable or
$\F=\Z_2$.
\\
(a)
If $M$ is embedded into a closed
3-manifold $Q$, then
$\dim H_1(Q;\F)\ge \dim H_1(M;\F)-~g$.
\\
(b) There
is a closed 3-manifold $Q$
containing $M$
such that $\dim H_1(Q;\F)= \dim H_1(M;\F)-~g$
and $Q$ is orientable
if $M$ is orientable.
\end{Th}

Part (a) is simple: it follows from the Mayer-Vietoris sequence, see the proof
at the end of the introduction.
Proof of part (b) (i.e., the construction of `minimal' $Q$)
is given in \S2. It is based on symplectic linear algebra
and Poincar\'e's theorem on the image of the mapping class group
of a surface $P$ in $\mathrm{Aut}(H_1(P;\Z))$. 

\begin{Th}
\label{Th2}
Let $M$ be a compact connected orientable
3-manifold with boundary.
Denote $g:=\rk H_1(\bd M;\Z)/2$.
\\
(a)
If $M$ is embedded into a closed
3-manifold $Q$, then
$H_1(Q;\Z)$ has
a sub-quotient isomorphic to
$$C(M):=\Z^{\rk H_1(M;\Z)-g}\oplus \mathrm{Tors}\ H_1(M;\Z).$$
(b) 
Suppose $H_1(M;\Z)\cong\Z^m$.
Then there
is a closed orientable 3-manifold $Q$
containing $M$
such that
$H_1(Q;\Z)\cong C(M)=\Z^{m-g}$.
\\
(c)
There is a compact connected orientable
3-manifold $M$ with boundary
which is not
embeddable into any
closed 3-manifold $Q$
such that
$H_1(Q;\Z)\cong C(M)$.
\end{Th}

Here $\rk X$ and $\Tor X$ are, respectively,
the rank and the torsion subgroup of an abelian group $X$.
Again, part (a) is essentially known
and part (b) is new; it is proved
after Theorem~\ref{Th1}(b) in \S2.
We present an example for part (c) in \S3.

{\it Remark.}
Suppose a closed orientable 3-manifold $Q$ contains $M$
and $H_1(Q;\Z)\cong C(M)$.
Then for {\it each} field $\F=\Z_p$
and $\F=\Q$ we get $\dim H_1(Q;\F)=\dim H_1(M;\F)-\rk H_1(\bd M;\Z)/2$,
while the proof of Theorem~\ref{Th1}(b) generally provides
{\it different} `minimal' manifolds
for different fields.

\smallskip

\begin{corollary}
\label{CorEmb}
Let $M$ be a compact orientable
3-manifold with boundary
and suppose
$G=\Z$, $G=\Z_p$ or $G=\Q$.
Then $M$ embeds into some $G$-homology
3-sphere if and only if
$H_1(M;G)\oplus H_1(M;G)\cong H_1(\bd M;G)$.
\end{corollary}

Corollary~\ref{CorEmb} is straightforward.
Corollaries~\ref{CorHomS} and~\ref{CorGrGen} are proved below in this section.
The construction of the `minimal' $Q$ in Corollary~\ref{CorGrGen}
is simpler than the general construction in Theorem~\ref{Th1}.
However, the lower estimation here is harder
and is reduced to the lower estimation in Theorem~\ref{Th1}
by the following lemma.
This lemma is proved in \S4.

\begin{lemma}
\label{LemGraph}
Let $L$ be a connected graph.
Suppose that the product $L\times S^1$ is embedded into
a 3-manifold $Q$.
Suppose that either
$Q$ is orientable or
$L$ is not homeomorphic to $S^1$ or $I$.
Then the regular neighborhood of $L\times S^1$ in $Q$
is homeomorphic to the product $K\times S^1$ for
a certain
2-manifold $K$ containing $L$.
If $Q$ is orientable, then $K$ is also orientable.
\end{lemma}

For instance, let $K_5$ be
the complete graph on 5 vertices.
Corollary~\ref{CorGrGen} implies that $K_5\times S^1$ is embeddable into a certain
closed orientable 3-manifold $Q$
such that $\dim H_1(Q;\F)=2$ and is not embeddable into any closed orientable
3-manifold
with the first homology group of dimension 0 or 1.
This result was obtained by A.~Kaibkhanov (unpublished).
The non-embeddability of $K_5\times S^1$ into $S^3$ was stated
by M.~Galecki and T.~Tucker (as far as the author knows, unpublished) and proved
by M.~Skopenkov in \cite{Sk03}.
\footnote
{
The non-embeddability
of $K_5\times S^1$ into $S^3$ could be proved in a simpler way
using the van Kampen theorem if we assumed
that $S^3\setminus U(K_5\times S^1)$
is homeomorphic to a disjoint union of
solid tori. (Here $U(K_5\times S^1)$ denotes the regular neighborhood
of $K_5\times S^1$ in $S^3$.)
However, this assumption is
not trivial to prove and becomes wrong if we replace
$K_5$ by some other graph $G$ such that
$G\times S^1$ embeds into $S^3$.
For example, let $G$ be a point.
Take a knotted embedding
$S^1\subset S^3$.
Then
$S^3\setminus U(S^1)$
is not homeomorphic to a solid torus.
}

\smallskip
The structure of the paper is as follows.
Now we prove Corollary~\ref{CorGrGen}, Theorems~\ref{Th1}(a) and~\ref{Th2}(a).
In this section we also prove Corollary~\ref{CorHomS},
for which we will need
Lemma~\ref{LemThick} below.
In \S2
we prove Theorems~\ref{Th1}(b) and~\ref{Th2}(b).
In \S3 we provide an example which proves Theorem~\ref{Th2}(c).
In \S4 we prove Lemmas~\ref{LemGraph} and~\ref{LemThick}.
The proof of both lemmas
uses the classification of 3-dimensional thickenings
of 2-polyhedra \cite{BRS99}.

\begin{example}
\label{ExampGr}
For $\F=\Z_p$ or $\F=\Q$
denote $r(M;\F):=\dim H_1(M;\F)- \dim H_1(\bd M;\F)/2$.
\\
(a) Let $\Xi$ be a surface of genus $g$ with $h$ holes.
Then $r(\Xi\times S^1;\F)=2g$.
\\
(b) Let $\Xi$ be a connected sum
of $k$ $\mathbb{R}P^2$'s
with $h$ holes.
Then $r(\Xi\times S^1;\Z_2)=k$.
\end{example}

\smallskip
{\it Proof of Corollary ~\ref{CorGrGen} modulo Theorem~\ref{Th1} and Lemma~\ref{LemGraph}.}
Since $S^1\times S^1\subset S^3$
and $I\times S^1\subset S^3$,
it suffices to consider the case when $L$
is not homeomorphic to $S^1$ or $I$.
Corollary~\ref{CorGrGen} now follows from 
Lemma~\ref{LemGraph},
Example~\ref{ExampGr}(a)
and
Theorem~\ref{Th1}.~$\blacksquare$

\smallskip
{\it Proof of Theorems \ref{Th1}(a) and \ref{Th2}(a).}
Suppose that $M\subset Q$,
where $M$ is a
3-manifold with boundary
and $Q$ is a closed 3-manifold.
In this paragraph, the 
homology coefficients are $\Z$, $\Z_p$ or $\Q$.
Let $i:H_1(\bd M)\to H_1(M)$,
$I:H_1(M)\to H_1(Q)$
be the inclusion-induced homomorphisms.
From the exact sequence of pair $(Q,M)$ 
we obtain that
$H_1(Q)$ has a subgroup isomorphic to
$H_1(M)/\Ker I$. 
From the Mayer-Vietoris sequence for $Q=M\cup_{\bd M}(Q-M)$
we obtain
$\Ker I \subset \img i$.
So
$H_1(M)/\img i$
is a quotient of
$H_1(M)/\Ker I$.

Let us prove Theorem~\ref{Th1}(a);
here the coefficients are $\F=\Z_p$
or $\F=\Q$.
By the known `half lives -- half dies' lemma,
$\dim \img i=\dim H_1(\bd M;\F)/2=g$, see
\cite[p.158]{FF89}, \cite[Lemma 3.5]{Ha}.
Thus $\dim H_1(Q;\F)\ge \dim H_1(M;\F)-g$.

To prove Theorem~\ref{Th2}(a),
it is left to check that
$C(M)\cong K:=H_1(M;\Z)/\img i$.
Indeed,
we obtain that
$\rk K=g$
by the universal coefficients formula
and the argument from the previous paragraph
for $\Q$-coefficients,
and
$\Tor K=\Tor H_1(M,\bd M;\Z)=\Tor H_1(M;\Z)$
by the exact sequence of pair
$(M,\bd M)$ and Poincar\'e duality.
$\blacksquare$

\smallskip
Let $P$ be a (finite) polyhedron.
If a 3-manifold $M$ is a regular neighborhood of $P\subset M$,
then the pair $(M,P)$ is called a {\it 3-thickening} of $P$
\cite{RoSa72}.
If we say that two thickenings are homeomorphic,
we mean that they are homeomorphic in the category of thickenings,
i.e. the homeomorphism in question
is relative to the polyhedron embedded into each thickening.

The following lemma
is known to specialists, but the author
has not found any proof
in literature.
This lemma is proved by combining \cite{BRS99} and \cite{Sk95}
(also see \cite{La00}); we prove it in \S4.

\begin{lemma}
\label{LemThick}
Each polyhedron $P$ has (up to homeomorphism) a finite number of
orientable 3-thickenings.
There exists an algorithm that
for a given polyhedron $P$
constructs
all its orientable 3-thickenings
(i.e., constructs their triangulations),
or tells that the polyhedron has none.
\end{lemma}

{\it Proof of Corollary \ref{CorHomS} modulo
Corollary~\ref{CorEmb} and Lemma \ref{LemThick}.}
Clearly, $P$ is embeddable into
an orientable 3-manifold $Q$ if and only if
there exists an
orientable 3-thickening
of $P$ which is embeddable into $Q$.
So the algorithm for Corollary~\ref{CorHomS} is as follows.
First, the algorithm constructs all
orientable 3-thickenings of $P$ with the help of Lemma~\ref{LemThick}.
If there are no such thickenings, then $P$ is not embeddable into any
orientable 3-manifold, and the algorithm gives the negative answer.
Otherwise, 
the algorithm
checks the condition of Corollary~\ref{CorEmb}
for each orientable 3-thickening of $P$
and gives the positive answer
if the condition was fulfilled for at least one 3-thickening.
$\blacksquare$

\smallskip
{\it Remark.}
Our methods do not lead to an algorithm
for embeddability of 2-polyhedra into $\R^3$
because we do not
deal with the fundamental group,
which is presumably much harder to do.
%

\section{Proof of Theorems~\ref{Th1}($b$),~\ref{Th2}($b$) (construction of a manifold $Q$)}

In this section give a proof of
Theorem~\ref{Th1}(b) and then
slightly modify it
to prove Theorem~\ref{Th1}(b).

\smallskip
{\it Proofs of Theorem \ref{Th1}(b).}
Denote $\F:=\Z_p$ or $\F:=\Q$.
In the current proof,
if coefficients in a
homology group are omitted,
they are assumed to be in $\F$.

Let $X\subset \R^3$ be the standardly embedded 
disjoint union of handlebodies such that
$\bd X\cong  \bd M$ 
and let $i:H_1(\bd M)\to H_1(M)$,
$i':H_1(\bd X)\to H_1(X)$
be the inclusion-induced homomorphisms.
We construct the required manifold $Q$
as a union of $X$ and $M$ along certain
diffeomorphism $f:\bd X\to \bd M$.
Consider the Mayer-Vietoris sequence
$$H_1(\bd M)\overset{i\oplus i' f_*^{-1}}\longrightarrow H_1(M)
\oplus H_1(X)\to H_1(Q)\to \widetilde{H_0}(\bd M)=0.$$
It follows that
$$H_1(Q)
\cong \frac{H_1(M)\oplus H_1(X)}{(i\oplus i' f_*^{-1})H_1(\bd M)}.$$
Suppose the map
$i\oplus i' f_*^{-1}$ 
is a monomorphism.
Then $\dim H_1(Q)=\dim H_1(M)-g$
as required.
So our goal now is to construct
a map $f:\bd X\to \bd M$ such that
$i\oplus i' f_*^{-1}$ is a monomorphism.

\smallskip
Let us introduce new notation.
For 
$G=\Z,\Q$ or $\Z_p$
a bilinear
form 
$\omega: G^{2g}\otimes G^{2g}\to G$ is called {\it symplectic}
if it is non-degenerate, skew-symmetric
and, when $G=\Z$, unimodular.
A submodule $B\subset G^{2g}$ will be called a ($G$-)\emph{Lagrangian} 
if $\omega|_B\equiv 0$ and $G^{2g}/B\cong G^g$.
Denote by $\Lin X$ the linear span of a subset $X$
of a vector space.
We will need the following lemma.

\begin{lemma}
Let $\omega$
be a symplectic form on $\Z^{2g}$.
\\
(a)
Denote by $\phi:\Z^{2g}\to{\Z_p}^{2g}$
the homomorphism
which applies $\mathrm{mod}\ p$
reduction to each component
and
by $\omega_{\Z_p}$
the symplectic form on ${\Z_p}^{2g}$
which is $\mathrm{mod}\ p$
reduction of $\omega$.
For each ${\Z_p}$-Lagrangian
$A\subset {\Z_p}^{2g}$
there exists
a $\Z$-Lagrangian $B\subset \Z^{2g}$
such that $\phi B=A$.
\\
(b)
Denote by $\phi:\Z^{2g}\to\Q^{2g}$
the inclusion
and by $\omega_\Q$
the symplectic form on $\Q^{2g}$
defined by the restriction $\omega_\Q|_{\Z^{2g}}\equiv \omega$.
For each $\Q$-Lagrangian
$A\subset \Q^{2g}$
there exists
a $\Z$-Lagrangian $B\subset \Z^{2g}$
such that $\Lin\phi B=A$.
\end{lemma}

{\it Proof of Lemma~2.1.}
Part (b) is obvious.
Let us prove part (a).
Recall that
$\phi$ is the reduction $\mathrm{mod}\ p$.
Take a set of generators
$\{e_i,f_i\}_{i=1}^g$
for $\Z^{2g}$
such that
$\omega(e_i,f_i)=\delta_{ij}$.
Then $\Lin\{\phi e_i\}_{i=1}^g$
is a $\Z_p$-Lagrangian,
hence there exists
$h_{\Z_p}\in\mathrm{Sp}(2g,{\Z_p})$
taking $\Lin\{\phi e_i\}_{i=1}^g$ to $A$
because
$\mathrm{Sp}(2g,{\Z_p})$
acts transitively on Lagrangians.
Since $\mathrm{mod}\ p$ reduction
maps $\mathrm{Sp}(2g,\Z)$
epimorphically onto
$\mathrm{Sp}(2g,{\Z_p})$
\cite[Theorem~VII.21]{Ne72},
we can find $h\in \mathrm{Sp}(2g,\Z)$
such that $\phi h=h_{\Z_p}$.
Then $B:=\{h e_i\}_{i=1}^{g}$
is the required $\Z$-Lagrangian.
$\blacksquare$

\smallskip
{\it Continuation of proof of Theorem~\ref{Th1}(b).}
Denote by
$\cap:H_1(\bd M;\Z)\times H_1(\bd M;\Z)\to\Z$
the intersection form
and by
$\cap_\F:H_1(\bd M)\times H_1(\bd M)\to\F$
the induced form
(as in Lemma~2.1);
$\cap|_\F$ coincides with
the $\F$-coefficients intersection form on 
$H_1(\bd M)$.
It is well known that 
$\dim \Ker i=g$ and
$\cap_\F|_{\Ker i}\equiv 0$
(the last assertion is analogous to
\cite[p.158]{FF89}).
In other words, $\Ker i$ is Lagrangian
with respect to $\cap_\F$.
Clearly,
there exists another Lagrangian
$A\subset H_1(\bd M)$
such that
$\Ker i\cap A=\{0\}$.
Let $\phi$ be the homomorphism from Lemma~2.1.
By Lemma~2.1(a)
or Lemma~2.1(b)
(depending on what coefficient field $\F$
we are working with)
we obtain
a Lagrangian submodule $B\subset H_1(\bd M;\Z)$
such that
$\Lin\phi B=A$
(if $\F=\Z_p$, this is equivalent to $\phi B=A$).
Notice that $\Lin\phi B=A$
implies that $\Ker i\cap\Lin\phi B=\{0\}$.

Recall the Poincar\'e theorem \cite{P} that
for a handlesphere $S$ every automorphism
of $H_1(S;\Z)$
preserving the intersection form
is induced by some self-diffeomorphism of $S$.

Denote $i'_\Z:\ H_1(\bd X;\Z)\to H_1(X;\Z)$
the inclusion-induced homomorphism; then
$\Ker i'_\Z$ is generated by the meridians
and is a $\Z$-Lagrangian
in $H_1(\bd X)$.
Thus there exists a diffeomorphism
$f:\ \bd X\to \bd M$
such that
$f_*\Ker i'_\Z=B$.
(Indeed, suppose that $\bd X\cong \bd M$ is connected.
Pick any
diffeomorphism $h_1:\bd X\to \bd M$.
Then $K:={h_1}_*\Ker i'_\Z\subset H_1(\bd M,\Z)$
is a $\Z$-Lagrangian.
By the Poincar\'e theorem
and because $\mathrm{Sp}(2g,\Z)$
acts transitively on $\Z$-Lagrangians
there exists
a self-diffeomorphism $h_2$ of $\bd M$
such that ${h_2}_*K=B$.
Now take $f:=h_2h_1$.
If $\bd M$ is not connected, apply this construction
componentwise.)

Because $X$ is a disjoint union of handlebodies,
$\Ker i'=\Lin \phi\Ker i'_\Z$
(if $\F=\Z_p$ and not $\Q$, then $\Ker i'=\phi\Ker i'_\Z$).
So
$$\Ker i' f_*^{-1}=f_* \Lin \phi\Ker i'_\Z=\Lin\phi f_*\Ker i'_\Z=\Lin\phi B.$$
Recall that $\Ker i\cap \Lin\phi B=\{0\}$.
Therefore
$i\oplus i' f_*^{-1}$ is monomorphic.
$\blacksquare$

\smallskip
{\it Proof of Theorem~\ref{Th2}(b).}
We use notation similar to the previous proof
and work with $\Z$-coefficients here.
Recall that $\Ker i$ is a $\Z$-Lagrangian,
i.e. $\cap|_{\Ker i}\equiv 0$ and $H_1(\bd M)/{\Ker i}\cong\Z^g$
\cite[p.158]{FF89};
thus we can find
a set of generators
$\{x_1,\ldots,x_{2g}\}\in H_1(\bd M)$
such that $\{x_1,\ldots,x_g\}$
generate $\Ker i$
and
$\{x_g,\ldots,x_{2g}\}$
also generate a Lagrangian.
Then there exists a diffeomorphism 
$f:\bd X\to \bd M$ such that
$\Ker i'f_*^{-1}$ is generated by
$\{x_{g+1},\ldots,x_{2g}\}$.
This is done analogously to
the proof of Theorem~\ref{Th1}(b)
using the Poincar\'e theorem
\footnote{This step is actually easier than
in Theorem~\ref{Th1}(b) because here we do not need Lemma~2.1}.
By construction we obtain
$$
H_1(Q)\cong\frac{H_1(M)\oplus H_1(X)}{(i\oplus i' f_*^{-1})H_1(\bd M)}
\cong
\frac{H_1(M)}{i H_1(\bd M)}\oplus
\frac{H_1(X)}{(i' f_*^{-1})H_1(\bd M)}\cong
\frac{H_1(M)}{i H_1(\bd M)}\cong C(M).
$$
The second
group in the direct sum is
obviously zero for
$X$ a disjoint union of handlebodies.
The last isomorphism is shown in the
proof of Theorems~\ref{Th1}(a),~\ref{Th2}(a).~$\blacksquare$



\section{Proof of Theorem~\ref{Th2}($c$)}

In this section we omit
$\Z$-coefficients. Theorem~\ref{Th2}(c)
is implied by the following two lemmas.

\begin{lemma}
\label{lemA}
There exists a connected orientable 3-manifold $M$ such that
\\
(1) $\bd M$ is a torus and
$H_1(M)\cong \Z\oplus\Z_2$.
\\
(2) 
Let $l$ and $m$ generate $H_1(M)$ and $2m=0$.
For some generators $a$, $b$
of $H_1(\bd M)\cong\Z\oplus\Z$ the inclusion-induced homomorphism
$i:H_1(\bd M)\to H_1(M)$ is given by
$i(a)=2l$, $i(b)=m$.
\end{lemma} 

\begin{figure}[t]
\centering
\includegraphics{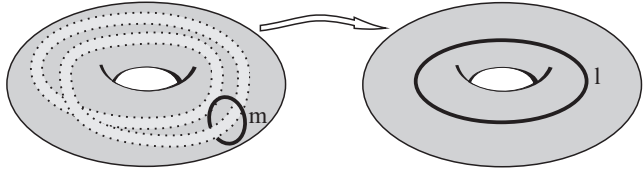}
\caption{Construction of a manifold in Lemma~\ref{lemA}.}
\label{Pict}
\end{figure}

\smallskip
{\it Proof.}
Let $D:=D^2\times S^1$ be a solid torus and $D'$ its copy.
Cut out from $D$ another solid torus which
lies inside $D$ and
runs twice along the parallel of $D$
(see Figure~1).
Glue the result
to $D'$ along $\bd D=\bd D'$.
It is easily seen that the orientable 3-manifold $M$ obtained satisfies
(1), (2).
The generators of $H_1(M)$ as in (2) are
shown on Figure~1.~$\blacksquare$

\begin{lemma}
\label{lemB}
Consider a manifold
$M$ from Lemma~\ref{lemA}. Then 
$C(M)=\Z_2$
(the group $C(M)$ is introduced in Theorem~\ref{Th2})
but
$M$ is not embeddable
into any closed 3-manifold $Q$
such that $H_1(Q)\cong~\Z_2$.
\end{lemma}

{\it Proof.}
Obviously, $C(M)=\Z_2$.
Suppose to the contrary that 
there is an embedding $M\subset Q$.
Denote by $X$ the closure of $Q\setminus M$
and by $i':H_1(\bd X)=H_1(\bd M)\to H_1(X)$ the
inclusion-induced homomorphism.
It follows from the Mayer-Vietoris sequence that
$$
H_1(Q)\cong \frac{H_1(M)\oplus H_1(X)}{(i\oplus i')H_1(\bd M)},
\quad\mbox{thus},
\quad
H_1(Q)\ \mbox{contains the subgroup}
\quad
R:=\frac{H_1(M)}{i(\Ker i')}.
$$
First, suppose $Q$ is orientable.
Then the rank of $\Ker i'$ is equal to 1,
so $\Ker i'$ is generated by $pa+qb$
for some $p,q\in \Z$.
Notice that $i(pa+qb)=2pl+qm$.
We obtain that $R$ is generated
by $l$ and $m$ with the following two relations:
$2m=0$, $2pl+qm=0.$
Clearly, $R\neq 0$
and $R\neq\Z_2$
since the determinant of the matrix
$
\left(\begin{smallmatrix}
0 & 2\\
2p & q
\end{smallmatrix}\right)
$
is divisible by $4$ but never equals $\pm 2$ or $\pm 1$,
as it should be when $R\cong\Z_2$ or $R=0$.

The case of non-orientable $Q$ is analogous.
We have now to consider the cases
$\rk\Ker i'=0$ and $\rk\Ker i'=2$.
In the first case, $R=\Z\oplus\Z_2$.
In the second case, the matrix of
relations for $R$:
$
\left( \begin{smallmatrix}
0 & 2p & 2r\\
2 & q & s
\end{smallmatrix}
\right)^t
$
is such that all of its 2$\times$2-minors 
are divisible by $4$.
This again implies that $R\neq 0$
and $R\neq \Z_2$.
$\blacksquare$

\smallskip
{\it Remark.}
The manifold $M$ constructed in Lemma~3.1 is embeddable
into a 3-manifold $Q$
with $H_1(Q)\cong \Z_2\oplus\Z_2$
and into
$S^1\times S^2$ with
$H_1(S^1\times S^2)\cong\Z$
(both manifolds are obtained by gluing a solid torus to $M$ appropriately).
These two manifolds verify Theorem~\ref{Th1}(b) for this particular manifold $M$:
the first manifold $Q$ when $\F\neq\Z_2$,
and $S^1\times S^2$ when $\F=\Z_2$.

\section{Proofs of Lemmas \ref{LemGraph}, \ref{LemThick}}

\smallskip
We will use results from \cite{BRS99};
let us state them here briefly and
prove Lemma~\ref{LemGraph} after that. The proof of
Lemma~\ref{LemThick} uses the same results and
is given at the end of this section.

\smallskip
{\it A classification of $3$-thickenings of 2-polyhedra \cite{BRS99}.}

Let $P$ be a 2-polyhedron. By $P'$
we will denote the 1-subpolyhedron which is the set of points in
$P$ having no neighborhood homeomorphic to 2-disk.
By $P''$ we will denote a (finite) set of points of $P'$
having no neighborhood homeomorphic to a book with $n$ sheets
for some $n\ge 1$. Take a point in any component of $P'$
containing no point of $P''$. Denote
by $F$ the union of $P''$ and these points.

Suppose that $\cup_{A\in F}\ \lk A$ is embeddable into $S^2$.
(Here $\lk$denotes link of a point.)
Take a collection of embeddings $\{g_A:\lk A\to S^2\}_{A\in F}.$
Take the closure $d\subset P'$
of a connected component of $P'\setminus P''$
and denote its ends by
$A,B\in F$ (possibly, $A=B$).
Now $d$ meets $\lk A\cup \lk B$ at two points
(distinct, even when $A=B$).
If for each such $d$ the maps $g_A$ and $g_B$
give the same or the opposite orders of rotation of the pages
of the book at $d$ then the collection $\{g_A\}$
is called {\it faithful}.
Two collections of embeddings
$\{f_A:\lk A\to S^2\},\ \{g_A:\lk A\to S^2\}$
are called {\it isopositioned}, if there
is a family of homeomorphisms $\{h_A:S^2\to S^2\}_{A\in F}$
such that $h_A\circ f_A=g_A$ for each $A\in F$.
This relation preserves faithfulness.
Denote by $E(P)$ the set of faithful collections up
to isoposition.

Suppose that $M$ is a 3-thickening of $P$.
Take any point $A\in F$ and consider its
regular neighborhood $R_M(A)$. Since $\bd R_M(A)$
is a sphere, we have a collection of embeddings
$\{\lk A\to\bd R_M(A)\}_{A\in F}$. Since for
each closure $d\subset P'$
of a connected component of $P'\setminus P''$
the regular neighborhood of $d$ is embedded into $M$,
this collection of embeddings is faithful. The
class $e(M)\in E(P)$ of this collection
is called the $e${\it-invariant}
of $M$. By $w_1(M)\in H^1(M;\Z_2)$ we denote the first Stiefel-Whitney class of $M$.

\begin{Th}
\cite[Theorem 3.1]{BRS99}.
Thickenings $M_1$, $M_2$ of $P$
are homeomorphic relative to $P$
if and only if
$w_1(M_1)|_P=w_1(M_2)|_P$ and
$e(M_1)=e(M_2)$.
\end{Th}

{\it Proof of Lemma~\ref{LemGraph}}
Without loss of generality we may assume that $Q$ is a regular
neighborhood of $L\times S^1$.
Due to Theorem 4.1, it is sufficient
to construct a 2-manifold $K$ containing
$L$ such that
\\
{\it(a)} $K\times S^1$ is a regular neighborhood
of $L\times S^1$, $e(K\times S^1)=e(Q)$ and
\\
{\it(b)} $w_1(K\times S^1)|_{L\times S^1}=w_1(Q)|_{L\times S^1}$.

First, let us construct a 2-manifold
$K$ satisfying {\it(a)}.
Take a triangulation of the graph $L$; we will work with this triangulation
only and denote it by the same letter $L$.
For each vertex $v$ in $L$ consider an arbitrarily oriented 2-disk $D^2_v$.
Consider the edges $e_1,\ldots,e_n$ containing $v$.
The embedding $L\times S^1\subset Q$ defines a cyclic
ordering of $e_1,\ldots,e_n$.
Take a disjoint union of $n$ arcs in $\partial D^2_v$
(each arc corresponding to an edge $e_i$) such that the cyclic ordering of the
arcs is the same as that of the edges.

For each edge $e$ connecting vertices
$u$ and $v$ connect $D^2_u$ and $D^2_v$ with a strip
$D^1\times D^1$,
gluing it along the two arcs that correspond to $e$.
The strip can be glued in two ways:
we can either twist it or not (with respect to
the orientations on $D^2_u$ and $D^2_v$).
After gluing a strip for each edge of $L$,
we get a union of disks and strips that
is a 2-manifold; denote it by $K$.
The manifold $K$ depends on choosing the twists.
However, any such $K$ satisfies {\it(a)},
no matter what the twists are.

By choosing the twists, let us
obtain the property {\it(b)}.

If $Q$ is orientable,
glue all the strips without twists. Then
$K$ is orientable, and
$w_1(K\times S^1)|_{L\times S^1}=w_1(Q)|_{L\times S^1}=0$.

Now let us choose
the twists in the other case:
$L$ is not homeomorphic
to $S^1$ or $I$ (and $Q$ is not necessarily orientable).
Denote the set of all edges of $L$ by $E$.
Take a point $O\in S^1$.
Take a set of cycles $c_1,\ldots,c_s\in Z_1(L;\Z_2)$ such that
$[c_1],\ldots,[c_s]\in H_1(L;\Z_2)$
is a basis.
Represent $w_1(Q)|_{L\times \{O\}}$
as a cochain $\{a_e\in\{0,1\}\}_{e\in E}$
so that
$\text{for all } k,\quad 1\le k\le s,\quad
\sum_{e\in c_k}a_e\ \mathrm{mod}\ 2=\langle w_1(Q)|_{L\times \{O\}},c_k\rangle$.
For each edge $e\in E$, twist
the corresponding strip
if $a_e=1$, and
do not twist the corresponding strip if $a_e=0$.
We now obtain $w_1(K\times S^1)|_{L\times \{O\}}=w_1(Q)|_{L\times \{O\}}$
by construction. We claim that the constructed $K$
satisfies {\it(b)}.

Indeed, take a vertex
$v$ of degree at least 3.
This can be done,
because $L$ is not homeomorphic to $S^1$ or $I$.
The homology classes of
$$c_i\times\{O\}, \quad 1\le i\le s,\quad\mbox{and}\quad \{v\}\times S^1$$
form a basis of $H_1(K\times S^1;\Z_2)$.
But
$$\langle w_1(Q),\{v\}\times S^1\rangle=0=
\langle w_1(K\times S^1),\{v\}\times S^1\rangle$$
because the regular neighborhood of $\{v\}\times S^1$ in $Q$
is orientable (the orientation is defined
by the orientation on $S^1$
and the cyclic ordering of the link of $v$ because $\mathrm{deg}\ v\ge 3$).
Thus we obtain $w_1(K\times S^1)|_{L\times S^1}=w_1(Q)|_{L\times S^1}$,
and the proof is finished.
$\blacksquare$

\smallskip

{\it Proof of Lemma \ref{LemThick}.}
Let $P$ be a 2-polyhedron.
We use the notation from the beginning of this section.
Take a faithful collection $\{g_A\}_{A\in F}$ of embeddings.
If the phrase from
the definition of faithfulness:
{\it
`the maps $g_A$ and $g_B$
give the same or the opposite orders of rotation of the pages
of the book at $d$'
}
is true even in the form
{\it
`the maps $g_A$ and $g_B$
always give
the opposite orders of rotation of the pages at $d$'},
then the collection $\{g_A\}$
is called {\it orientably faithful}.
Two collections
$\{f_A\},\{g_A\}$
are called {\it orientably isopositioned}, if there
is a family of
orientation-preserving
homeomorphisms $\{h_A:S^2\to S^2\}_{A\in F}$
such that $h_A\circ f_A=g_A$ for each $A\in F$.
This relation preserves
the property of being orientably faithful.
Denote by $SE(P)$ the set of
orientably faithful collections up
to orientable isoposition.

An orientable 3-thickening $M$ of $P$
induces an \emph{$se$-invariant} $se(M)\in SE(P)$.
It is an oriented version
of the $e$-invariant and is defined analogously.
The following is essentially proved in \cite{Sk95} and \cite{La00}:
every class
$c\in SE(P)$
is an $se$-invariant of
some orientable 3-thickening of $P$.
These papers give an algorithm for
construction of
such thickening.
Moreover, if two orientable 3-thickenings $M_1,M_2$ of $P$
have the same $se$-invariants $se(M_1)=se(M_2)\in SE(P)$, they are
homeomorphic
(this follows from Theorem~3,
since the Stiefel-Whitney classes are zeros in the orientable case).

The set $SE(P)$ is obviously finite.
Hence the number of orientable 3-thickenings of $P$ is finite.
The algorithm for construction of all
orientable 3-thickenings of $P$
is as follows. For each class $c\in SE(P)$
build a corresponding orientable
$3$-thickening using the
construction from \cite{Sk95}, \cite{La00}. Theorem~4.1
guarantees that we will obtain all orientable
3-thickenings as result.
$\blacksquare$


\bigskip
{\bf Acknowledgements.}
The author is grateful to
D.~Crowley, S.~Melikhov 
and especially to A.~Skopenkov for useful discussions.

\end{document}